\documentclass[12pt,reqno]{amsart}

\usepackage{txfonts}
\usepackage{amsmath,amssymb,amsthm} 
\usepackage{hyperref}
\usepackage[numbers]{natbib}
\usepackage{rotating}
\usepackage{adjustbox}
\usepackage{array}      
\usepackage{multicol}         
\usepackage[bottom]{footmisc}
\usepackage{float}
\usepackage{cancel}
\usepackage[dvipsnames]{pstricks}
\usepackage{pst-grad} 
\usepackage{pst-plot}  
\usepackage{multirow} 
\usepackage{color}
\usepackage{setspace}
\usepackage{colortbl}
\usepackage{wrapfig}
\usepackage{lmodern}
\usepackage{tcolorbox}
\usepackage[pangram]{blindtext}
\usepackage{graphicx,epstopdf}
\usepackage{parskip}
\usepackage[left=2.8cm, right=2.8cm, top=3.0cm, bottom=3.0cm]{geometry}
\usepackage{natbib}



\newcommand{\Tt}{\mathcal{T}}

\newtheorem{theorem}{Theorem}[section]
\newtheorem{lemma}[theorem]{Lemma}
\newtheorem{proposition}[theorem]{Proposition}
\newtheorem{corollary}[theorem]{Corollary}
\theoremstyle{definition}
\newtheorem{definition}[theorem]{Definition}

\newtheoremstyle{boldtitle-plainbody}
  {\topsep}
  {\topsep}
  {\normalfont}
  {0pt}
  {\bfseries}
  {.}
  { }
  {}
\theoremstyle{boldtitle-plainbody}
\newtheorem{remark}[theorem]{Remark}
\newtheorem{example}[theorem]{Example}

\title[Minimum number of non-monochromatic simplices]{On the minimum number of non-monochromatic simplices for Sperner labelings of a regular triangulation}

\author[L. Á. CALVO PASCUAL]{LUIS ÁNGEL CALVO}
\address{Luis Ángel Calvo Pascual\\
Dpto. de Métodos Cuantitativos, ICADE\\
Universidad Pontificia Comillas de Madrid, Spain}
\email{lacalvo@comillas.edu}

\author[S. MERCHÁN]{SUSANA MERCHÁN}
\address{Susana Merchán\\
Dpto. de Matemática e Informática Aplicadas a las Ingenierías Civil y Naval\\
Universidad Politécnica de Madrid, Spain}
\email{susana.merchan@upm.es}

\author[D. RABOSO]{DULCINEA RABOSO}
\address{Dulcinea Raboso\\
Dpto. de Matemáticas y Ciencias de Datos\\
Universidad San Pablo CEU de Madrid, Spain}
\email{dulcinea.raboso@gmail.com}

\author[J. RODRIGO]{JAVIER RODRIGO}
\address{Javier Rodrigo\\
Dpto. de Matemática Aplicada\\
Universidad Pontificia Comillas de Madrid, Spain}
\email{jrodrigo@comillas.edu}

\author[J. S. RODRÍGUEZ]{JOSÉ SAMUEL RODRÍGUEZ}
\address{José Samuel Rodríguez\\
Dpto. de Matemáticas, Instituto Gran Capitán, Madrid, Spain}
\email{jose.rodriguez219@educa.madrid.org}

\date{\today}

\begin{document}

\begin{center}
\begin{minipage}{0.97\textwidth}
\begin{abstract}
\small Attending to an open problem in the literature stated by Mirzakhani and Vondrák, we give a lower bound of the number of non-monochromatic simplices for Sperner labelings of the vertices of a triangulation of a given $ k$-simplex with vertices of integer coordinates. This triangulation maximizes the number of simplices over all the triangulations of the $ k$-simplex with vertices of integer coordinates.
\end{abstract}
\vspace{0.5cm}

\maketitle

\vspace{0.3cm}
\end{minipage}
\end{center}

\vspace{0.5cm}
\section{Introduction}	

\noindent The search for lower bounds of the number of non-monochromatic simplices (at least two labels) in Sperner labelings of the vertices of a triangulation of a simplex has deserved attention by its connection with the hypergraph labeling problem.

This problem aims to find a cut in a hypergraph that minimizes the sum of the assignment costs of the vertices and the weight of the edges. It is related to computational complexity problems such as uniform geometric labeling or monotone restricted MSCA (see \cite{Ene}, \cite{Chek}, \cite{Klein} for more details).   
\cite{Mirza2} (see also \cite{Mirza}) succeeds in the problem of finding a tight lower bound on the number of non-monochromatic simplices in Sperner labelings of the vertices of a complete hypergraph contained in a $ k$-simplex. This hypergraph was not a simplicial subdivision of the $ k$-simplex.

An open problem discussed in \cite{Mirza2} is finding a lower bound for the number of non-monochromatic simplices in Sperner labelings of simplicial subdivisions (triangulations) of the simplex \( \Delta_{k,q} \), where  \( \Delta_{k,q} \) denotes the set of integer vectors in \( \mathbb{Z}^k \) with non-negative entries that sum to $q$. Specifically, on the final page of \cite{Mirza2}, Mirzakhani and Vondrák proposed the following conjecture:
\begin{quote}
\textit{``For a Sperner-admissible labeling of a regular simplicial subdivision [...], what is the minimum possible number of non-monochromatic cells? [...]. We conjecture that for a fixed \( j \leq k \) and as \( q \to \infty \), the number of cells containing at least \( j \) colors is on the order of \(\Omega( q^{k-j}) \).''}
\end{quote}

\noindent
In this paper, we address this problem by working with the concept of \emph{regular triangulation}, $\Tt$,  a symmetric subdivision of $\Delta_{k,q}$ into simplices of equal volume, whose vertices are the integer lattice points.

We establish the following main results:

\begin{theorem} \label{prop42}
The minimum number of non-monochromatic simplices for a Sperner labeling of the vertices of the regular triangulation \( \Tt \) of \( \Delta_{k,q} \) satisfies
\[
m_{k,q} \geq \binom{q + k - 3}{k - 2}.
\]
\end{theorem}

 This bound is proved to be tight for the initial cases $q=1$ or $k=2$. For $q=2$ we give a better lower bound with an improvement of the multiplicative constant.
Next, we give an upper bound for the minimum number of non-monochromatic simplices, which is $O(q^{k-2})$, with the multiplicative constant depending on the dimension $k$. Concretely,

\begin{theorem} \label{prop51}
The minimum number of non-monochromatic simplices for a Sperner labeling of the vertices of the regular triangulation \( \Tt \) of \( \Delta_{k,q} \) satisfies
\[
m_{k,q} \leq q^{k-1} - (q - 1)^{k-1}.
\]
\end{theorem}

We achieved the bounds for the triangulation of $\Delta_{k,q}$ with the vertices set $\Delta_{k,q}\cap \mathbb{Z}^k$ defined in \cite{Mirza2}.

This paper is organised as follows. Section \ref{preliminaries} introduces definitions and notation used throughout. In Section \ref{Section 3}, we present and characterise the regular triangulation $\Tt$. Section \ref{Section 4} establishes a lower bound on the minimum number of non-monochromatic simplices within this triangulation. In Section \ref{Section 5}, we determine an upper bound. Finally, in Section \ref{Section 6}, we give some conclusions and future directions of work.

\section{Basic concepts and notation}\label{preliminaries}

In this section, we fix notation and recall some definitions.

\textbf{Notation}
\begin{itemize}
    \item 
    $\Delta_k$ is a simplex of $k$ vertices, $k\geq 2$.
    \item
    $Conv\{v^1,\ldots, v^k\}$ is the convex hull of the vertices $v^1,\ldots, v^k\,.$
    \item
    $x_i$ is the $i$-th coordinate of the vector $x\in\mathbb{Z}^k. $
\end{itemize}

\vspace{0.5cm}

\begin{definition} 
For $k, q\in\mathbb{N}$, $k\geq 2$, we define the $k$\emph{-simplex} $\Delta_{k,q}$ as
 $$\Delta_{k,q} = \left\{ (x_1, \ldots, x_k) \in \mathbb{Z}^k \,\middle|\, 
\sum_{i=1}^k x_i = q,\; x_i \geq 0, \text{ for all } i = 1, \ldots, k \right\}$$
$$=Conv\{(q,0,\ldots, 0), \ldots, (0,\ldots,0, q)\}.$$
\end{definition}
\vspace{0.5cm}

\begin{definition} \label{SimplicialSubdivision}
A \emph{simplicial subdivision} (or \emph{triangulation}) $T_r$ of a $k$-simplex $\Delta_k$ is a finite collection of simplices (called \emph{cells}) satisfying the following conditions:
\begin{itemize}
    \item The union of all simplices in $T_r$ is exactly $\Delta_k$, i.e., 
    \[
    \bigcup_{\Delta \in T_r} \Delta = \Delta_k.
    \]
    \item For any pair of simplices $\Delta^1, \Delta^2 \in T_r$, the intersection $\Delta^1 \cap \Delta^2$ is either empty or a common face of both simplices.
\end{itemize}
\end{definition}
\vspace{0.5cm}

We give the definition of Sperner labeling as stated in \cite{Mirza}. Despite it is not as general as the classical one, it gives a good fit to our development.
\vspace{0.5cm}

\begin{definition} \label{SpernerLabeling}
Let $V$ denote the set of vertices of a triangulation of $\Delta_{k,q}$. A labeling $c: V \to \{1, \ldots, k\}$ is called a \emph{Sperner labeling} if, for every vertex $v = (v_1, v_2, \ldots, v_k) \in V$, the following condition holds
\[
v_i = 0 \quad \Rightarrow \quad c(v) \neq i.
\]
\end{definition}
\vspace{0.2cm}

Now, we define a particular case of Sperner labeling (see \cite{Mirza2}).
\vspace{0.5cm}

\begin{definition} \label{FirstChoiceLabeling}
The \emph{first choice labeling} is a Sperner labeling defined on the set of vertices $V$ as follows: for each vertex $v = (v_1, v_2, \ldots, v_k) \in V$, 
$$
c(v) = \min\left\{ i \in \{1, \ldots, k\} \mid v_i > 0 \right\}.
$$
\end{definition}
\vspace{0.5cm}

\section{The regular triangulation}\label{Section 3}

We will employ the simplicial subdivision of $\Delta_{k,q}$ defined in~\cite{Mirza2} (see also~\cite{Doug,Edel}), which we refer to as the \textit{regular triangulation}~$\Tt$. Before proceeding, we introduce the following preliminary definitions.
\vspace{0.5cm}

\begin{definition}
Let \( k, q \in \mathbb{N} \) with \( k \geq 2 \). The simplex \( R_{k,q} \) is given by
\[
R_{k,q} = \left\{ y \in \mathbb{R}^{k-1} \,\middle|\, 0 \leq y_1 \leq y_2 \leq \cdots \leq y_{k-1} \leq q \right\}.
\]
\end{definition}

We denote by \( W_{k,q} \) the set of integer lattice points contained in \( R_{k,q} \), that is,
$$
W_{k,q} = R_{k,q} \cap \mathbb{Z}^{k-1}.
$$
\vspace{0.2cm}

\begin{definition}
Let $w = (w_1, \ldots, w_{k-1}) \in W_{k,q}$. A permutation $\pi: \{1, \ldots, k-1\} \to \{1, \ldots, k-1\}$ is said to be \emph{consistent with} $w$ if 
\[
w_i = w_{i+1} \quad \Rightarrow \quad \pi(i) < \pi(i+1).
\]
\end{definition}
\vspace{0.5cm}

\begin{definition}\label{def3.3}
Let $\pi: \{1, \ldots, k-1\} \to \{1, \ldots, k-1\}$ be a permutation consistent with $w \in W_{k,q}$. We define the simplex
\[
\sigma(w, \pi) = \left\{ y \in \mathbb{R}^{k-1} \,\middle|\, y_i \geq 0, \quad 0 \leq (y - w)_{\pi(1)} \leq \cdots \leq (y - w)_{\pi(k-1)} \leq 1 \right\}.
\]
\end{definition}
\vspace{0.5cm}

\begin{definition} \label{triangulation_Tk}
The \emph{triangulation} $T_k$ of the region $R_{k,q}$ is the collection of all simplices $\sigma(w, \pi)$ arising from pairs $(w, \pi),$ where $\pi$ is a permutation consistent with $w \in W_{k,q}$, i.e.,
\[
T_k = \left\{ \sigma(w, \pi) \,\middle|\, w \in W_{k,q},\; \pi \text{ consistent with } w \right\}.
\]
\end{definition}

\vspace{0.5cm}

\begin{definition}\label{regular triangulation}
The \emph{regular triangulation} $\Tt$ of the discrete simplex $\Delta_{k,q}$ is the collection of simplices
\[
\Tt = \left\{ \phi(\sigma(w, \pi)) \,\middle|\, w \in W_{k,q},\; \pi \text{ consistent with } w \right\},
\]
where $\phi: R_{k,q} \to \Delta_{k,q}$ is the mapping given by
\[
\phi(y_1, \ldots, y_{k-1}) = \left( y_1,\; y_2 - y_1,\; y_3 - y_2,\; \ldots,\; y_{k-1} - y_{k-2},\; q - y_{k-1} \right).
\]
\end{definition}
\vspace{0.5cm}

To characterise the regular triangulation $\Tt$, we first analyse $T_k$ through the introduction of an associated graph.
\vspace{0.5cm}

\begin{definition}
The graph $G_k'$ is the pair $(V_k', E_k')$, where:
\begin{itemize}
    \item The vertex set is $V_k' = W_{k,q}$.
    \item The edge set $E_k'$ consists of pairs $\{w^1, w^2\} \subset W_{k,q}$ such that there exist $(w, \pi)$ that satisfy $w^1, w^2 \in \sigma(w, \pi).$
\end{itemize}
\end{definition}
\vspace{0.5cm}

\begin{definition}\label{defT_k'}
We define the collection $T_k'$ as the set of convex hulls of $k$ vertices in $W_{k,q}$ that are pairwise adjacent in the previous graph, that is,
\[
T_k' = \left\{ Conv(v^1, \ldots, v^k) \,\middle|\, v^i \in W_{k,q},\; \text{and } \{v^i, v^j\} \in E_k' \text{ for all } 1 \leq i < j \leq k \right\}.
\]
\end{definition}
\vspace{0.5cm}

We now establish the relationship between the sets $T_k$ and $T_k'$ introduced in Definitions~\ref{triangulation_Tk} and~\ref{defT_k'}, respectively.
\vspace{0.5cm}

\begin{proposition}\label{prop31}
It is satisfied that $T_k=T_k'$.
\end{proposition}

\begin{proof}
Let $\sigma(w, \pi) \in T_k$ be a simplex. By construction, its $k$ vertices $v^1, \ldots, v^k$ belong to $W_{k,q}$ and lie within $\sigma(w, \pi)$, therefore $v^1,\ldots,v^k$ are pairwise adjacent in $G'$ and then
$$Conv(v^1,\ldots,v^k)=\sigma(w,\pi)\in T'_k,$$
which shows that $T_k \subset T_k'$.

Conversely, suppose there exist $k$ pairwise adjacent vertices $v^1,\ldots,v^k$ such that 
$$Conv(v^1,\ldots,v^k)\notin T'_k.$$
Since the union of all simplices in $T_k$ covers $R_{k,q}$, i.e.,
\[
R_{k,q} = \bigcup_{(w, \pi)} \sigma(w, \pi),
\]
the convex hull $Conv(v^1, \ldots, v^k)$ must intersect the interior of some simplex $\sigma(w, \pi) \in T_k$. Since this edge is an intersection of simplices of $T_k$, then $\sigma(w',\pi)$ contains an interior point of $\sigma(w,\pi)$ for some $w'\in W_{k,q}$, a contradiction because the intersection of simplices of a triangulation is a common face. Therefore, our assumption must be false, and it follows that $T_k' \subset T_k$ and then $T_k=T'_k$ as desired.
\end{proof}

The previous Proposition \ref{prop31} implies that $T'_k$ is a triangulation of $R_{k,q}$. Now we give a manageable definition of $T'_k$.
\vspace{0.5cm}

\begin{definition}\label{def:Gk''}
The graph $G_k''$ is the pair $(V_k'', E_k''),$ where
\begin{itemize}
    \item The vertex set is $V_k'' = W_{k,q}$.
    \item For distinct vertices $w^1, w^2 \in W_{k,q}$, the edge $\{w^1, w^2\}$ belongs to $E_k''$ if either
    \[
    (w^1 - w^2)_i \in \{0, 1\}, \quad \text{for all } i = 1, \ldots, k-1,
    \]
    or
    \[
    (w^1 - w^2)_i \in \{-1, 0\}, \quad \text{for all } i = 1, \ldots, k-1.
    \]
\end{itemize}
\end{definition}
\vspace{0.5cm}

\begin{definition}
The collection $T_k''$ is the set of convex hulls of $k$ vertices in $W_{k,q}$ that are pairwise adjacent in the graph $G_k''$. That is,
\[
T_k'' = \left\{ Conv(v^1, \ldots, v^k) \,\middle|\, v^i \in W_{k,q},\; \{v^i, v^j\} \in E_k'', \text{ for all } 1 \leq i < j \leq k \right\}.
\]
\end{definition}
\vspace{0.5cm}

\begin{proposition}\label{310}
It is satisfied that $G_k'=G_k''$.
\end{proposition}

\begin{proof}
Since $V'_k=V''_k=W_{k,q}$, it suffices to prove that $E'_k=E''_k$. Assume that $\{w^1,w^2\}\in E'_k$, then $w^1,w^2\in\sigma(w,\pi),$ for some $w\in W_{k,q}$ and then
$$0\leq (w^1-w)_{\pi(1)}\leq\ldots\leq(w^1-w)_{\pi(k-1)}\leq 1.$$ 
Since $(w^1-w)_{\pi(i)}\in\mathbb{Z},$ for all $i=1,\ldots, k-1$, it follows that there exist and index $i$ such that:
$$(w^1-w)_{\pi(1)}=\ldots=(w^1-w)_{\pi(i)}=0, (w^1-w)_{\pi(i+1)}=\ldots=(w^1-w)_{\pi(k-1)}=1.$$

Similarly, for \( w^2 \), there exists an index \( j \) such that 
$$(w^2-w)_{\pi(1)}=\ldots=(w^2-w)_{\pi(j)}=0, (w^2-w)_{\pi(j+1)}=\ldots=(w^2-w)_{\pi(k-1)}=1.$$
Assume without loss of generality that $j>i$, then
\begin{itemize}
    \item for $s\leq i$, $$(w^1-w^2)_{\pi(s)}=(w^1-w)_{\pi(s)}-(w^2-w)_{\pi(s)}=0,$$

\item for $i<s\leq j$, $$(w^1-w^2)_{\pi(s)}=(w^1-w)_{\pi(s)}-(w^2-w)_{\pi(s)}=1-0=1,$$
\item for $j<s\leq k-1$, $$(w^1-w^2)_{\pi(s)}=(w^1-w)_{\pi(s)}-(w^2-w)_{\pi(s)}=1-1=0.$$ \end{itemize}
Therefore, $\{w^1,w^2\}\in E''_k.$

Now we show that $E''_k\subset E'_k$ by induction on $k$. For $k=2$,  if $\{w^1,w^2\}\in E''_2$, then we can assume without loss of generality that $w^2-w^1=1$. The permutation $\pi(1)=1$ is trivially admissible and satisfies
$$0\leq(w^1-w^2)_{\pi(1)}\leq 1,$$
so $w^1, w^2\in\sigma(w^1,\pi)$ and then $\{w^1,w^2\}\in E'_2$.  

Assume that $E''_{k-1}\subset E'_{k-1}$ and consider $G''_k=(V''_k,E''_k)$ and $\{w^1,w^2\}\in E''_k$. 
We can assume without loss of generality that $(w^1-w^2)_i\in\{0,1\}$ for $i=1,\ldots,k-1$.  We have two cases.

\textbf{\emph{Case 1.}} There exists an index $i$ such that $(w^1)_i=(w^1)_{i+1}$ and $(w^2)_i=(w^2)_{i+1}$.

We consider
$$w'^1=(w^1_1,\ldots,w^1_{i}, w^1_{i+2},\ldots,w^1_{k-1}), \hspace{1cm} w'^2=(w^2_1,\ldots,w^2_{i}, w^2_{i+2},\ldots,w^2_{k-1})\in\mathbb{Z}^{k-2}\,.$$
These satisfy that 
$$0\leq(w'^1)_1=(w^1)_1\leq(w'^1)_2\leq\ldots\leq(w'^1)_{k-2}\leq q,$$ 
$$0\leq(w'^2)_1=(w^2)_1\leq(w'^2)_2\leq\ldots\leq(w'^2)_{k-2}\leq q,$$ so $w'^1,w'^2\in R_{k-1,q}\cap\mathbb{Z}^{k-2}=W_{k-1,q}$. 

Since \( (w'^1 - w'^2)_i \) corresponds to $(w^1 - w^2)_j$, for some \( j \in \{1, \ldots, k-2\} \), and \( (w^1 - w^2)_j \in \{0,1\} \) by assumption, it follows that \( (w'^1 - w'^2)_i \in \{0,1\}\), for all \( i = 1, \ldots, k-2 \). Consequently, \( \{w'^1, w'^2\} \in E''_{k-1} \), and by the induction hypothesis, \( \{w'^1, w'^2\} \in E'_{k-1} \). Therefore, there exists \( w' \in W_{k-1,q} \) and a permutation \( \pi \) such that \( w'^1, w'^2 \in \sigma(w', \pi) \). We also have that $E''_{k-1}=E'_{k-1}$ and then $G''_{k-1}=G'_{k-1}$, so $T''_{k-1}=T'_{k-1}$ and the simplices of the triangulation of $R_{k-1,q}$ are the convex hulls of $k-1$ vertices pairwise adjacent in $G''_{k-1}$.

The map
$$f(x_1,\ldots,x_i,x_{i+1},\ldots,x_{k-2})=(x_1,\ldots,x_i,x_i,x_{i+1},\ldots,x_{k-2})$$ defines a bijection  between $W_{k-1,q}$ and  $W_{k,q}\cap\{(x_1,\ldots,x_{k-1}) \,\mid\, x_i=x_{i+1}\}$, where the inverse is given by
$$f^{-1}(x_1,\ldots,x_i,x_i,x_{i+1},\ldots,x_{k-2})=(x_1,\ldots,x_i,x_{i+1},\ldots,x_{k-2}).$$
Since $f$ just duplicates a coordinate and $f^{-1}$ removed the repeated coordinate, it follows that edges of $E''_{k-1}$ are mapped to edges of $E''_{k}$ whose vertices lies in $$W_{k,q}\cap\{(x_1,\ldots,x_{k-1})\mid x_i=x_{i+1}\}.$$ Conversely, edges of $E''_{k}$ with vertices in $W_{k,q}\cap\{(x_1,\ldots,x_{k-1})\mid x_i=x_{i+1}\}$ are mapped back in edges of $E''_{k-1}$  by $f^{-1}$. 
Therefore, $G''_{k-1}$ is isomorphic to the induced subgraph $$G''_{k}\big|_{W_{k,q}\cap\{(x_1,\ldots,x_{k-1})\mid x_i=x_{i+1}\}}.$$

In this way, the set of convex hulls of $k-2$ vertices of $W_{k,q}\cap\{(x_1,\ldots,x_{k-1})\mid x_i=x_{i+1}\}$ pairwise adjacent in $G''_{k}$ is a triangulation of $R_{k,q}\cap\{(x_1,\ldots,x_{k-1})\mid x_i=x_{i+1}\}$, we call it $T$, and $w^1,w^2$ share a simplex of $T$, since $w'^1,w'^2$ share a simplex of $T''_{k-1}$.

Moreover, for each \( \sigma(w, \pi) \in T'_k \) containing \( k-1 \) vertices of
$$
W_{k,q} \cap \{ (x_1, \ldots, x_{k-1}) \in \mathbb{Z}^{k-1} \mid x_i = x_{i+1} \},
$$
the intersection
$$
\sigma(w, \pi) \cap \{ (x_1, \ldots, x_{k-1}) \in \mathbb{R}^{k-1} \mid x_i = x_{i+1} \}
$$
is a simplex of
$$
R_{k,q} \cap \{ (x_1, \ldots, x_{k-1}) \in \mathbb{R}^{k-1} \mid x_i = x_{i+1} \}.
$$
Since these $k-1$ vertices are pairwise adjacent in $G''_k$, it follows that
$$
\sigma(w, \pi) \cap \{ (x_1, \ldots, x_{k-1}) \in \mathbb{R}^{k-1} \mid x_i = x_{i+1} \} \in T.
$$
We denote these simplices by $\sigma'(w, \pi)$.
  
We have 
$$\bigcup_{(w,\pi)}\{\sigma'(w,\pi)\cap\mathbb\{(x_1,\ldots,x_{k})\mid x_i=x_{i+1}\}=R_{k,q}\cap\mathbb\{(x_1,\ldots,x_{k})\mid x_i=x_{i+1}\}$$

Otherwise, 
$$R_{k,q}\cap\mathbb\{(x_1,\ldots,x_{k}) \mid x_i=x_{i+1}\}-\bigcup_{(w,\pi)}\{\sigma'(w,\pi)\cap\mathbb\{(x_1,\ldots,x_{k})\mid x_i=x_{i+1}\}$$ 
is a nonempty open subset of $R_{k,q}\cap\mathbb\{(x_1,\ldots,x_{k}) \mid x_i=x_{i+1}\}$, so it has dimension $k-2$. 
But since 
$$\bigcup_{(w,\pi)}\{\sigma(w,\pi)\cap\mathbb\{(x_1,\ldots,x_{k}) \mid x_i=x_{i+1}\}=R_{k,q}\cap\mathbb\{(x_1,\ldots,x_{k}) \mid x_i=x_{i+1}\},$$
we would have that 
$$R_{k,q}\cap\mathbb\{(x_1,\ldots,x_{k}) \mid x_i=x_{i+1}\}-\bigcup_{(w,\pi)}\{\sigma'(w,\pi)\cap\mathbb\{(x_1,\ldots,x_{k})\mid x_i=x_{i+1}\}$$ is the union of a finite number of sets with dimension $<k-2$, a contradiction.

We also have that 
$$(\sigma'(u,\pi)\cap\mathbb\{(x_1,\ldots,x_{k})\mid x_i=x_{i+1}\})\cap(\sigma'(v,\pi)\cap\mathbb\{(x_1,\ldots,x_{k})\mid x_i=x_{i+1}\})$$ is either empty or a common face for $u,v\in W_{k,q}$ , $u \neq v$  since $\sigma'(u,\pi)\cap\sigma'(v,\pi)$ is either empty or a common face for $u,v\in W_{k,q}$, $u \neq v$.

It follows that $\{\sigma'(w,\pi)\cap\mathbb\{(x_1,\ldots,x_{k})\mid x_i=x_{i+1}\}\}$ is a triangulation of 
$$R_{k,q}\cap\mathbb\{(x_1,\ldots,x_{k})\mid x_i=x_{i+1}\}$$
  included in $T$, and then 
  $$T=\left\{\sigma'(w,\pi)\cap\mathbb\{(x_1,\ldots,x_{k})\mid x_i=x_{i+1}\}\right\}.$$
  
This implies that $w^1,w^2\in\sigma'(w,\pi),$  for some $w\in W_{k,q}$ , and then $\{w^1,w^2\}\in E'_k$. Therefore,  $E''_k\subset E'_k$ in this case.

\textbf{\emph{Case 2.}} Suppose there does not exist an index $i$ such that  $$(w^1)_i=(w^1)_{i+1}\hspace{1cm}\text{and}\hspace{1cm}(w^2)_i=(w^2)_{i+1}.$$

We can assume that $i_1<i_2<\ldots <i_j$ , $i_{j+1}<i_{j+2}<\ldots <i_{k-1,}$  for some $j$ and $(w^1-w^2)_{i_s}=0$   for $1\leq s\leq j$, $(w^1-w^2)_{i_s}=1$  for $j<s\leq k-1$. 

We define the following permutation of $\{1,\ldots ,k-1\}$: $\pi (t)=i_{t}$ . Let us see that $\pi (t)$  is consistent with  $w^2$ if either $(w^2)_j \neq (w^2)_{j+1}$  or $(w^2)_j=(w^2)_{j+1}$  and $i_{j}<i_{j+1}$.

Indeed, $\pi (t)$ is strictly increasing except for (maybe) $j,j+1$, so if  $(w^2)_j \neq (w^2)_{j+1}$, then $\pi (t)$ is consistent with $w^2$. If  $(w^2)_j = (w^2)_{j+1}$ and $i_j<i_{j+1}$  (that is to say, if $(w^1)_i = (w^2)_i$  for $i\leq j$), then  $\pi (j)<\pi (j+1)$  and $\pi (t)$ is strictly increasing so $\pi (t)$ is consistent with $w^2$.
 
Moreover, we have
$$(w^1-w^2)_{\pi(1)}=\ldots =(w^1-w^2)_{\pi(j)}=0,$$
 $$(w^1-w^2)_{\pi(j+1)}=\ldots =(w^1-w^2)_{\pi(k-1)}=1,$$ it follows that $w^1\in \sigma(w^2,\pi)$ and then $\{w^1,w^2\}\in E_k'$ in this case.
 
If $(w^2)_{j}=(w^2)_{j+1}$  and $i_{j+1}<i_{j}$, then $(w^1)_j \neq (w^1)_{j+1}$ by the assumption.

Since $(w^1-w^2)_t\in \{0,1\}$ for every $t$, the former implies that $(w^1)_j = (w^2)_{j}$  and $(w^1)_{j+1} = (w^2)_{j+1}+1$ , so $i_1\leq j<j+1\leq i_k$ , and then the permutation 
$$\pi(1)=i_2,\ldots ,\pi(j-1)=i_j,\pi(j)=i_1,\pi(j+1)=i_k,\pi(j+2)=i_{j+1},\ldots ,\pi(k)=i_{k-1}$$ 
satisfies that if $(w^2)_{i}=(w^2)_{i+1},$  then $\pi(i)<\pi(i+1),$ except for the cases $i=j-1$, $i=j+1$. So, if $(w^2)_{j-1} \neq (w^2)_{j}$  and $(w^2)_{j+1} \neq (w^2)_{j+2}$, then $\pi$ is admissible.

In the case $(w^2)_{j+1} = (w^2)_{j+2}$, we obtain
$$(w^1)_{j+2} \geq (w^1)_{j+1}=(w^2)_{j+1}+1=(w^2)_{j+2}+1.$$
As a result, $(w^1)_{j+2}= (w^2)_{j+2}+1=(w^2)_{j+1}+1=(w^1)_{j+1}$, a contradiction.

Similarly, if $(w^2)_{j-1} = (w^2)_{j}$, then 
$$(w^1)_{j}= (w^2)_{j}=(w^2)_{j-1}\leq (w^1)_{j-1},$$ so $(w^1)_{j} = (w^1)_{j-1},$ a contradiction. Therefore $\pi$ is admissible, with 
$$(w^1-w^2)_{\pi(1)}=\ldots =(w^1-w^2)_{\pi(j)}=0,$$
$$(w^1-w^2)_{\pi(j+1)}=\ldots =(w^1-w^2)_{\pi(k-1)}=1,$$ it follows that $w^1\in \sigma(w^2,\pi)$ and then $\{w^1,w^2\}\in E_k'$, so $E''_{k}\subset E'_{k}$ as desired.

\end{proof}
\vspace{0.5cm}

\begin{corollary}\label{rem1}
The triangulations \( T_k \), \( T_k' \), and \( T_k'' \) are equal.
\end{corollary}

\begin{proof}
By Proposition~\ref{310}, the graphs \( G_k' \) and \( G_k'' \) are isomorphic, then the sets of \( k \) vertices that are pairwise adjacent in \( G_k' \) and \( G_k'' \) coincide. Thus, the collections of simplices defined as the convex hulls of such sets are equal, that is, \( T_k' = T_k'' \).

On the other hand, by Proposition~\ref{prop31}, we have \( T_k = T_k' \).

Combining these equalities, it follows that
\[
T_k = T_k' = T_k''.
\]
\end{proof}

\begin{definition}\label{def:GK}
The graph $G_k $ is the pair $(V_k, E_k)$, where
\begin{itemize}
    \item The vertex set is $V_k = \Delta_{k,q} \cap \mathbb{Z}^k$.
    \item For distinct vertices $v^1, v^2 \in V_k$, the edge $\{v^1, v^2\}$ belongs to $E_k$ if 
    \[
    (v^1 - v^2)_i \in \{-1, 0, 1\} \quad \text{for all } i = 1, \ldots, k,
    \]
    and the number of entries equal to $1$ and $-1$ are equal and when reading the nonzero entries from left to right, the signs alternate.
\end{itemize}
\end{definition}

We now establish the relationship between the graphs $G''_k$ and $G_k$ introduced in Definitions~\ref{def:Gk''} and~\ref{def:GK}, respectively.

\begin{proposition}\label{isom}
$G_k''$ is isomorphic to $G_k$.  
\end{proposition}
\begin{proof}
The mapping $\phi|_{W_{k,q}}: W_{k,q}\rightarrow\Delta_{k,q}\cap Z^k$ is biyective, with inverse 
$$\phi^{-1}(x_1,\ldots ,x_k)=(x_1,x_1+x_2,x_1+x_2+x_3,\ldots ,x_1+x_2+\ldots +x_{k-1}).$$ 
If \( \{w^1, w^2\} \in E_k'' \), we may assume without loss of generality that \( (w^1 - w^2)_i \in \{0,1\} \) for all \( i = 1, \ldots, k-1 \). Let \( i_0 \) denote the smallest index such that \( (w^1 - w^2)_{i_0} = 1 \).

If \( i_0 > 1 \), then for all \( i < i_0 \) we have
$$
(\phi(w^1))_i = (w^1)_i - (w^1)_{i-1} = (w^2)_i - (w^2)_{i-1} = (\phi(w^2))_i,
$$
and
$$
(\phi(w^1))_{i_0} = (w^1)_{i_0} - (w^1)_{i_0 - 1} = (w^2)_{i_0} + 1 - (w^2)_{i_0 - 1} = (\phi(w^2))_{i_0} + 1.
$$
If \( i_0 = 1 \), the relation simplifies to
$$
(\phi(w^1))_{i_0} = (w^1)_{i_0} = (w^2)_{i_0} + 1 = (\phi(w^2))_{i_0} + 1.
$$

Next, let \( i_1 \) be the smallest index greater than \( i_0 \) such that \( (w^1 - w^2)_{i_1} = 0 \), if such an index exists. If \( i_1 > i_0 + 1 \), then for all \( i_0 < i < i_1 \) we have
$$
(\phi(w^1))_i = (w^1)_i - (w^1)_{i-1} = (w^2)_i + 1 - \left( (w^2)_{i-1} + 1 \right) = (\phi(w^2))_i,
$$
and at \( i_1 \) we have
$$
(\phi(w^1))_{i_1} = (w^1)_{i_1} - (w^1)_{i_1 - 1} = (w^2)_{i_1} - \left( (w^2)_{i_1 - 1} + 1 \right) = (\phi(w^2))_{i_1} - 1.
$$

  Following the same reasoning with $i_2$, $i_3, \ldots$ we conclude that the entries with $1$'s and $-1$'s in $\phi(w_1)-\phi(w_2)$ alternate, and we have the same number of $1$'s and $-1$'s since
  $$(\phi(w_1)-\phi(w_2))_1+\ldots +(\phi(w_1)-\phi(w_2))_k=q-q=0.$$ 
  
This implies that $\{\phi(w^1),\phi(w^2)\}\in E_k$.
 
Conversely, if $\{v^1,v^2\}\in E_k$, then we can assume without loss of generality that the lowest index $i$ such that $(v^1-v^2)_i\neq 0$, $i_1$, satisfies $(v^1-v^2)_{i_1}=1$.
If $i_j$, $j=1, ..., 2t$, is the ordered list of indices such that $(v^1-v^2)_{i_j}\neq 0$, then we have $(v^1-v^2)_{i_j}=1,$ for $j$  an odd number and $(v^1-v^2)_{i_j}=-1$ for $j$ an even number.
Consequently, if $i_1>1$, then, for $i<i_1$:
$$(\phi^{-1}(v^1))_i=\sum_{j=1}^i(v^1)_j=\sum_{j=1}^i(v^2)_j=(\phi^{-1}(v^2))_i,$$ where
$$(\phi^{-1}(v^1))_{i_1}=\sum_{j=1}^{i_1}(v^1)_j=\sum_{j=1}^{i_1-1}(v^1)_j+(v^1)_{i_1}=\sum_{j=1}^{i_1-1}(v^2)_j+(v^2)_{i_1}+1=(\phi^{-1}(v^2))_{i_1}+1.$$
 
Moreover, for the indices $i$ such that $i_1<i<i_2$  (if any), we obtain
$$(\phi^{-1}(v^1))_i=\sum_{j=1}^i(v^1)_j=\sum_{j=1}^{i_1}(v^1)_j+\sum_{j={i_1+1}}^i(v^1)_j=\sum_{j=1}^{i_1}(v^2)_j+1+\sum_{j={i_1+1}}^i(v^2)_j=(\phi^{-1}(v^2))_i+1,$$
with 
$$(\phi^{-1}(v^1))_{i_2}=\sum_{j=1}^{i_2}(v^1)_j=\sum_{j=1}^{i_2-1}(v^1)_j+(v^1)_{i_2}=\sum_{j=1}^{i_2-1}(v^2)_j+1+(v^2)_{i_2}-1=(\phi^{-1}(v^2))_{i_2}.$$

Iterating the process with $i_3$, $i_4$, $\ldots$, we obtain that $(\phi^{-1}(v^1)-\phi^{-1}(v^2))_i\in \{0,1\}$ for $i=1,...,k-1$, and then $\{\phi^{-1}(v^1),\phi^{-1}(v^2))\} \in E''_k$.

So $G_k''$ is isomorphic to $G_k$ as desired.
\end{proof}

\begin{corollary}
The regular triangulation $\Tt$ of $\Delta_{k,q}$ can be described as
\[
\Tt = \left\{ Conv(v^1, \ldots, v^k) \,\middle|\, v^1, \ldots, v^k \in \Delta_{k,q} \cap \mathbb{Z}^k,\; \text{and } \{v^i, v^j\} \in E_k, \text{ for all } 1 \leq i < j \leq k \right\}.
\]
\end{corollary}

\begin{proof}
By Definition \ref{regular triangulation}, $$\Tt=\{\phi(\sigma(w,\pi)) \mid \sigma(w,\pi)\text{ is a simplex of } T_k\}.$$ According to Corollary \ref{rem1}, the simplices of $T_k$ are $Conv(w^1,\ldots, w^k),$  for $w^1,\ldots, w^k$ pairwise adjacent in $G_k''$. Since $\phi$ defines a graph isomorphism between $G_k''$ and $G_k$ (as stated in Proposition~\ref{isom}), then the image of each such simplex under $\phi$ corresponds to a convex hull of $k$ vertices in $\Delta_{k,q},$ that are pairwise adjacent in $G_k$. It follows that the regular triangulation $\Tt$ consists precisely of the simplices $Conv(v^1, \ldots, v^k)$ for vertices $v^1, \ldots, v^k$ in $\Delta_{k,q} \cap \mathbb{Z}^k$ that are pairwise adjacent in $G_k$, which proves the result.
\end{proof}

\section{The lower bound}\label{Section 4}
In order to establish a lower bound $\Omega(q^{k-2})$ for the number of non-monochromatic simplices, we need a preliminary definition of a hypergraph included in \cite{Mirza}.
\vspace{0.3cm}

\begin{definition}\label{simplex-lattice}
The \emph{Simplex-Lattice Hypergraph} $H_{k,q}$ is a pair $ (V_{k,q}, E_{k,q})$, where:
\begin{itemize}
    \item The vertex set is $V_{k,q} = \Delta_{k,q} \cap \mathbb{Z}^k$.
    \item The hyperedge set is
    \[
    E_{k,q} = \left\{ \left\{ (b_1+1, b_2, \ldots, b_k),\, (b_1, b_2+1, \ldots, b_k),\, \ldots,\, (b_1, b_2, \ldots, b_k+1) \right\} \,\middle|\, (b_1, \ldots, b_k) \in V_{k,q-1} \right\}.
    \]
\end{itemize}
\end{definition}
\vspace{0.2cm}

\begin{definition}
Let $T$ be a triangulation of $\Delta_{k,q},$ whose vertices are labeled with a Sperner labeling \( c: V \to \{1, \ldots, k\} \). A simplex in \( T \) is said to be \emph{monochromatic} if all its vertices are assigned the same label by $c$; otherwise, the simplex is called \emph{non-monochromatic}.
\end{definition}

We denote by \( m_{k,q} \) the \emph{minimum number of non-monochromatic simplices} that appear in any Sperner labeling of the regular triangulation $\Tt$ of \( \Delta_{k,q} \).
\vspace{0.3cm}

\begin{theorem}[\textbf{Theorem \ref{prop42} of the Introduction}]
The number of non-monochromatic simplices for a Sperner labeling of the vertices of the regular triangulation \( \Tt \) of \( \Delta_{k,q} \) satisfies
\[
m_{k,q} \geq \binom{q + k - 3}{k - 2}.
\]
\end{theorem}

\begin{proof}
The simplices associated with the hyperedges in $E_{k,q}$ are simplices in the regular triangulation $\Tt$. Indeed, observe that for any two vertices 
\[
(b_1, \ldots, b_i+1, \ldots, b_j, \ldots, b_k) \quad \text{and} \quad (b_1, \ldots, b_i, \ldots, b_j+1, \ldots, b_k),
\]
the difference is
\[
(0, \ldots, 0, -1, 0, \ldots, 0, 1, 0, \ldots, 0),
\]
which corresponds to a vector of the form described in the definition of the graph $G_k$ (see Definition \ref{def:GK}). Then, all vertices in each hyperedge of $E_{k,q}$ form a simplex in $\Tt$.

Therefore, the number of non-monochromatic simplices in $\Tt$ is at least the number of non-monochromatic hyperedges in $E_{k,q}$. According to Proposition 2.1 of~\cite{Mirza}, this number is bounded below by $\binom{q + k - 3}{k - 2}$.
\end{proof}

Let us see an example for which the lower bound of Theorem \ref{prop42} is not tight.
\vspace{0.5cm}

\begin{example}\label{ex41} For $q = 2,$ any Sperner labeling of the vertices of the regular triangulation $\Tt$ has at most one monochromatic simplex.

Indeed, if every monochromatic simplex has label $k$ (for instance), then the monochromatic simplices have the vertices either in $x_k = 1$ or in $x_k = 2$. But for $q = 2,$ there is just one simplex with this property: the simplex with $k - 1$ vertices in $x_k = 1$ and a vertex in $x_k = 2$, namely $(0, \ldots, 0, 2)$.

If there exist monochromatic simplices with labels $1$ and $k$ (for instance), then we would have only one simplex with label $1$ and only one simplex with label $k$ as we have seen before, so the $k-1$ vertices in  $V_k\cap\{x_1=1\}$ would have label $1$ and the $k-1$ vertices in $V_k\cap\{x_k=1\}$ would have label $k$, so we obtain vertices in 
$$V_k\cap\{(x_1\ldots,x_k)|\,\,x_1=x_k=1\}$$
with two labels, a contradiction.

Therefore, the number of non-monochromatic simplices is at least $2^{k-1}-1$, and then $m_{k,2}\geq 2^{k-1}-1$.
This bound improves the lower bound of Theorem \ref{prop42}, for $k>2$.
\end{example}

\section{An upper bound of $m_{k,q}$} \label{Section 5}
Now we see an upper bound of $m_{k,q}$.
\begin{theorem}[\textbf{Theorem \ref{prop51} of the Introduction}]
The minimum number of non-monochromatic simplices for a Sperner labeling of the vertices of the regular triangulation \( \Tt \) of \( \Delta_{k,q} \) satisfies
\[
m_{k,q} \leq q^{k-1} - (q - 1)^{k-1}.
\]
\end{theorem}
\begin{proof}
Consider the first choice labeling (Definition \ref{FirstChoiceLabeling}). In this labeling, a simplex is non-monochromatic if and only if it contains at least one vertex with first coordinate \( x_1 = 0 \). These are all the simplices of the regular triangulation $\Tt$ except for the simplices such that their vertices have first coordinate $x_1\geq 1$. 

Since the simplices of the regular triangulation $\Tt$ with a vertex in $x_1=0$ are included in $x_1\leq 1$, then the simplices of $\Tt$ such that their vertices have first coordinate $x_1\geq 1$ triangulate:
$$\{(x_1,\ldots, x_k)|x_1\geq 1, x_2\geq 0,\ldots, x_k\geq 0, x_1+\ldots +x_k=q\}\sim\Delta_{k,q-1},$$
so the number of said simplices is $(q-1)^{k-1}$ and then the number of non-monochromatic simplices for the first choice labeling is $q^{k-1}-(q-1)^{k-1}$. This implies that $$m_{k,q}\leq q^{k-1}-(q-1)^{k-1},$$ as desired.
\end{proof}

 \begin{remark}
Theorem \ref{prop51} implies that 
$$m_{k,q}\leq q^{k-1}-(q-1)^{k-1}\sim (k-1)q^{k-2}.$$
Note that the multiplicative constant of this upper bound tends to infinity as $k$ tends to infinity. In addition to this, the multiplicative constant of the lower bound of Theorem \ref{prop42} tends to $0$ as $k$ tends to infinity.
 \end{remark}
\vspace{0.5cm}

\begin{remark}
For $k=2$, $q=1$ or $q=2$, the lower bounds of Theorem \ref{prop42} and Example \ref{ex41} meet the upper bound of Theorem \ref{prop51} and we get $$m_{2,q}=m_{k,1}=1;\hspace{1cm}  m_{k,2}=2^{k-1}-1.$$
\end{remark}

\section{Conclusions and future lines}\label{Section 6}

We establish a lower bound in the open problem  of determining the minimum number ($m_{k,q}$) of non-monochromatic simplices for Sperner labelings of the vertices of a simplicial subdivision of $\Delta_{k,q}$ given in \cite{Mirza}.

\begin{figure}[H] 
    \centering
    \includegraphics[width=0.8\textwidth]{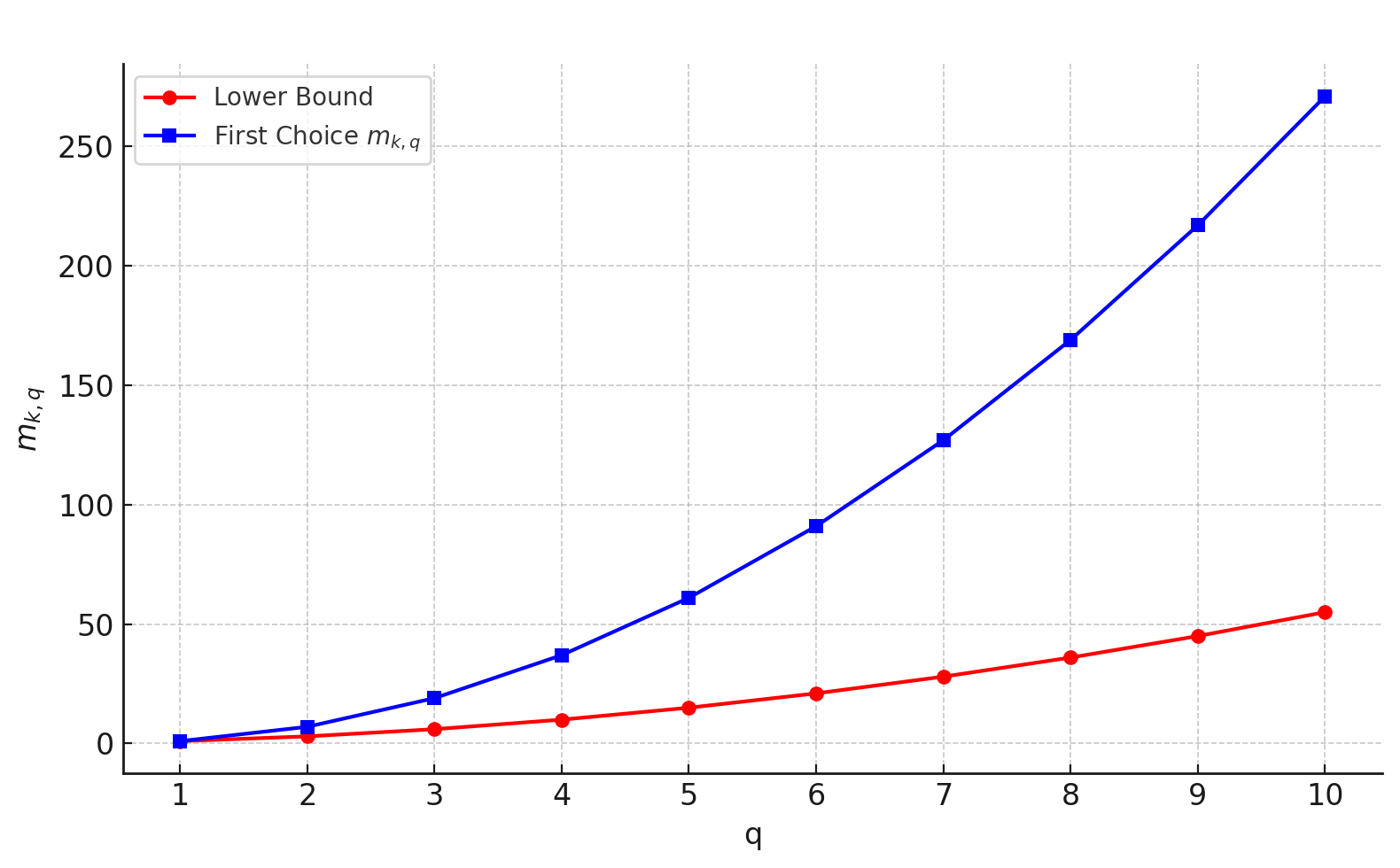}
    \caption{Comparison of Lower Bound and First Choice $m_{k,q},$ for $k=4$.}
    \label{figure1}
\end{figure}

In that paper, the authors stated that it was crucial for their application to find a lower bound of $m_{k,q}$  with a tight multiplicative constant. A natural direction for future work could be to narrow the gap between the multiplicative constant $\frac{1}{(k-2)!}$ of the lower bound of $m_{k,q}$ (Theorem \ref{prop42}) and the multiplicative constant $k-1$ of the upper bound of $m_{k,q}$ (Theorem \ref{prop51}). See Figure~\ref{figure1} for a visualization of this gap.

In the initial cases for which we have established the exact value of $m_{k,q}$ ($k=2$, $q=1$, $q=2$), this value coincides with the upper bound of Proposition \ref{prop51} and is attained for the first choice labeling.  We conjecture that this is the general case.

Our results suggest several research directions for fair division, hypergraph coloring, and combinatorial topology.

Determining the exact constant, or even the lower-order terms, in $m_{k,q}$ would sharpen worst-case guarantees for discrete fair-division algorithms. Classical rental-harmony models based on Sperner’s lemma~\cite{Su1999}, bounded envy-free cake-cutting protocols~\cite{aziz2016discrete}, and recent studies on mixed-resource division~\cite{liu2024mixed} all rely on understanding the structure of worst-case instances. Identifying the true minimum number of non-monochromatic simplices would help characterise these challenging cases and support the design of more robust allocation methods.

 Refining the value of $m_{k,q}$ could impact rainbow generalisations of the KKM lemma and multi-cut problems such as necklace splitting and the Hobby-rice theorem, where multiple balanced partitions are sought simultaneously. Key developments in this direction include degree-theoretic and combinatorial proofs of Sperner-type results~\cite{LeVan1982,Atanassov1996,kumar17,article}, polytopal generalisations~\cite{DELOERA20021}, purely combinatorial approaches~\cite{MEUNIER20061462}, rainbow and criticality extensions~\cite{Asada2018,Kaiser2024}, and recent homotopy-based methods~\cite{Dulinski2020}.

Finally, the regular triangulation considered here is closely related to the edgewise subdivision of a simplex~\cite{Edel} and to earlier work on simplicial mesh generation~\cite{Doug}. A more detailed geometric analysis of these constructions could lead to new algorithmic insights and improved combinatorial bounds.

\begin{thebibliography}{99}

\bibitem{Ene} A. Ene, J. Vondrák (2014) {\it Hardness of Submodular Cost Allocation: Lattice Matching and a Simplex Coloring Conjecture}, In Proc. of APPROX, 144--159.

\bibitem{Chek} C. Chekuri, A. Ene (2011) {\it Submodular Cost Allocation Problem and Applications}, Proc. of ICALP, 354--366.

\bibitem{Klein} J. M. Kleinberg, E. Tardos (2002) {\it Approximation Algorithms for Classification Problems with Pairwise Relationships: Metric Labeling and Markov Random Fields}, Journal of the ACM, 49, 5, 616--639.

\bibitem{Mirza2} M. Mirzakhani, J. Vondrák (2017) {\it Sperner’s Colorings and Optimal Partitioning of the Simplex}, A Journey Through Discrete Mathematics, Springer, 615--631.

\bibitem{Mirza} M. Mirzakhani, J. Vondrák (2015) {\it Sperner’s Colorings, Hypergraph Labeling Problems and Fair Division}, Proc. of ACM-SIAM SODA, 873--886.

\bibitem{Doug} W. M. Douglas (1992) {\it Simplicial Mesh Generation with Applications}, Ph.D. thesis, Cornell University.

\bibitem{Edel} H. Edelsbrunner, D. R. Grayson (2000) {\it Edgewise Subdivision of a Simplex}, Discrete \& Computational Geometry, 24, 4, 707--719.

\bibitem{article} T. Le, C. Van, N.-S. Pham, C. Saglam (2022) {\it A Direct Proof of the Gale–Nikaido–Debreu Lemma Using Sperner’s Lemma}, Journal of Optimization Theory and Applications, 194.

\bibitem{LeVan1982} C. Le Van (1982) {\it Topological Degree and the Sperner Lemma}, Journal of Optimization Theory and Applications, 37, 3, 371--377.

\bibitem{kumar17} A. Kumar, K. G. (2017) {\it The Simplex Reminiscent of Sperner's Lemma}, IARJSET, 4.

\bibitem{Su1999} F. E. Su (1999) {\it Rental Harmony: Sperner's Lemma in Fair Division}, The American Mathematical Monthly, 106, 10, 930--942.

\bibitem{MEUNIER20061462} F. Meunier (2006) {\it Sperner labellings: A combinatorial approach}, Journal of Combinatorial Theory, Series A, 113, 7, 1462--1475.

\bibitem{DELOERA20021} J. A. De Loera, E. Peterson, F. E. Su (2002) {\it A Polytopal Generalization of Sperner's Lemma}, Journal of Combinatorial Theory, Series A, 100, 1, 1--26.

\bibitem{aziz2016discrete} H. Aziz, S. Mackenzie (2016) {\it A Discrete and Bounded Envy-Free Cake Cutting Protocol for Any Number of Agents}, Proceedings of the 57th Annual IEEE Symposium on Foundations of Computer Science (FOCS), 416--427.

\bibitem{liu2024mixed} S. Liu, X. Lu, M. Suzuki, T. Walsh (2024) {\it Mixed Fair Division: A Survey}, Journal of Artificial Intelligence Research, 80, 1373--1406.

\bibitem{Atanassov1996} K. T. Atanassov (1996) {\it On Sperner's Lemma}, Studia Scientiarum Mathematicarum Hungarica, 32, 1, 71--74.

\bibitem{Asada2018} M. Asada, F. Frick, V. Pisharody, M. Polevy, D. Stoner, L. H. Tsang, Z. Wellner (2018) {\it Fair Division and Generalizations of Sperner- and KKM-type Results}, SIAM Journal on Discrete Mathematics, 32, 1, 591--610.

\bibitem{Dulinski2020} W. Duliński (2020) {\it Homotopies and Transcendental Extensions in Colouring Problems}, available at \url{https://arxiv.org/abs/2011.12273}.

\bibitem{Kaiser2024} T. Kaiser, M. Stehlík, R. Škrekovski (2024) {\it Criticality in Sperner's Lemma}, Combinatorica, 44, 5, 1041--1051.

\end{thebibliography}
\end{document}